\documentclass[11pt, leqno]{article}
\title{\bf Relations of the spaces \\$A^p(\Omega)$ and $C^p(\partial\Omega)$}
\author{V. Mastrantonis}
\date{}
\usepackage[greek,english]{babel}  
\usepackage[utf8x]{inputenc}
\usepackage{hyperref}
\RequirePackage{latexsym} \RequirePackage{amsthm}
\RequirePackage{amsmath}
\newcommand{\norm}[1]{\left\lVert#1\right\rVert}
\usepackage{enumerate}
\RequirePackage{amssymb} \RequirePackage{makeidx}
\usepackage{dsfont}
\usepackage{mathrsfs}
\usepackage[margin= 3.5cm]{geometry}

\newtheorem{theorem}{Theorem}[section]
\newtheorem{lemma}[theorem]{Lemma}
\newtheorem{proposition}[theorem]{Proposition}

\newtheorem{corollary}[theorem]{Corollary}
\newtheorem{definition}[theorem]{Definition}

\newtheorem{remark}[theorem]{Remark}

\numberwithin{equation}{section}

\begin{document}

\maketitle


\begin{abstract}
\footnotesize
In this paper we prove that for functions $f\in A(D)$ there is an equivalence between the continuous extension of their derivatives over the boundary and the differentiability  of the map $t\mapsto f(e^{it})$. More specifically, we are able to prove that $A^p(D)= A(D)\cap C^p(\mathbb{T})$ by making use  of the Poisson representation. Moreover, we extend our results over Jordan domains bounded by an analytic Jordan curve by using what we initially prove on the disk in combination with the Osgood- Caratheodory theorem.
\end{abstract}

\section{\textbf{Introduction}}
For every function $f\in A(D)$, where $D$ stands for the open unit disk in $\mathbb{C}$ and $A(\Omega)$ denotes the collection of functions defined on the Jordan Domain $\Omega$ that are holomorphic on $\Omega$ and continuous on $\overline{\Omega}$, we know that $f'$ is well defined in $D$. But what one say if the function $f'$ continuously extends on $\overline{D}$? Does this imply that the map $t \mapsto f(e^{it})$ is also differentiable and what is the relation between its derivative and the extension of $f'$ over the unit circle? Conversely, if the map $t\mapsto f'(e^{it})$ is differentiable what can be said about the derivative $f'$ in the open unit disk? Is it continuously extendable over $\overline{D}$? What happens if we replace the disk with a Jordan domain bounded by an analytic Jordan curve? The present paper deals with these types of questions.

We define, for any Jordan domain $\Omega$, the spaces $A^p(\Omega)$, $p\in \{0, 1, 2,... \}\cup\{\infty\}$ as the class of all holomorphic in $\Omega$ functions  $f$ whose derivatives $f^{(l)}$, $l\in \{0, 1, ...\}, l\leq p$, have continuous extensions over $\overline{\Omega}$. In section 2 we prove that for any $f\in A(D)$, the fact that $f\in A^p(D)$ is equivalent to the fact that the map $t \mapsto f(e^{it})$ is in $C^p(\mathbb{T})$, where $\mathbb{T}$ denotes the unit circle and as a corollary we have that $A^p(D)= A(D)\cap C^p(\mathbb{T})$. To prove this we use the Poisson representation. Furthermore, we are able to prove analogous results for functions $f\in A(D)$ whose derivatives continuously extend on an open arc of the unit circle and not on the whole circle. Moreover, by making use of Riemann's mapping theorem we are able to extend this result to the class of functions $f$ which are holomorphic in $D$, not necessarily in $A(D)$, but continuously extendable over an open arc of the unit circle, but once again not necessarily on the whole boundary. More specifically, for such an $f$ one can show that the first $p$ derivatives of the function continuously extend over the open arc $J$ if and only if the map $t\mapsto f(e^{it})$ belongs in $C^p(I)$, where $I= (a, b)$ is an interval such that $J= \{e^{it}: a< t< b \}$.

In the last section, we generalize the previous results over Jordan domains bounded by a Jordan curve, part of which is an analytic curve. More specifically, if $f$ is a holomorphic function defined in a Jordan domain $\Omega$ and $\gamma: I \mapsto\mathbb{C}$ is an analytic curve, $I= (a, b)\subset\mathbb{R}$ with $\gamma(I)\subset\partial\Omega$, then once again the first $p$ derivatives of $f$ continuously extend over $\gamma^{\ast}= \gamma(I)$ if and only if the map $t\mapsto f(\gamma(t))$ is in $C^p(I)$. In order to prove this we apply the Osgood-Caratheodory theorem \cite{[4]} and we make use of the previous results on the disk. If $\gamma(I)= \partial\Omega$, we use any conformal parametrization. In particular, the parametrization induced by Riemann mapping \cite{[1]} from $D$ onto the interior of $\gamma^\ast$ is a proper parametrization and can be used as well. The previous results can be extended to the case of domains of finite connectivity bounded by a finite set of disjoint Jordan curves, provided a segment of the boundary is an analytic arc.

\section{The case of the Disk.}
For $0 \leq p\leq +\infty$, we denote by $A^p(D)$ the space of holomorphic functions on $D$ whose derivatives of order $l\in\mathbb{N}, 0\leq l\leq p$, extend continuously over $\overline{D}$. It is topologized via  the semi-norms:
\[
|f|_l= \sup_{z\in\overline{D}}|f^{(l)}(z)|= \sup_{z\in\mathbb{T}}|f^{(l)}(z)|, 0\leq l \leq p, l\in\mathbb{N}
\]
For any $z= re^{i\theta}\in\mathbb{C}$ we denote by $P_z(t)$ or $P_r(\theta)$ the Poisson kernel \cite{[1]}.
\\ \newline
The above spaces can be defined for open sets $\Omega$ more general than the disk $D$. However, the simplest case is the case of the disk. Consequently, we firstly try to investigate the relations between functions in the classes $A^p(D)$ and $C^p(\mathbb{T})$, where $\mathbb{T}= \partial D$ is the unit circle. For that we will need the following:
\begin{lemma}
Let $u:\mathbb{T}\to \mathbb{C}$ such that $u\in C^p(\mathbb{T})$ and define:
\begin{equation}
U(z)= \frac{1}{2\pi} \int_{-\pi}^{\pi}u(t)P_z(t)dt
\end{equation}
then U is well defined in $\mathbb{C}\setminus\mathbb{T}$ and $C^\infty$ harmonic. The following holds for all $\theta_0\in\mathbb{R}$ and $0<r<1$:
\begin{equation}
\frac{d^lU}{d\theta^l}(re^{i\theta_0})= \frac{1}{2\pi}\int_{-\pi}^{\pi}u^{(l)}(t)P_r(\theta_0- t) dt
\end{equation}
for $l=0, 1, ...,p$ and hence:
\begin{equation}
\lim_{r\to1^-} \frac{d^lU}{d\theta^l}(re^{i\theta_0})= u^{(l)}(\theta_0)
\end{equation}
uniformly for all $\theta_0\in\mathbb{R}$.
\begin{proof}
We later prove, in Proposition 3.1, that a more general class of functions defined in a similar way are well defined and $C^{\infty}$ harmonic in $\mathbb{C}\setminus\mathbb{T}$. Notice that $U(re^{i\theta})= \frac{1}{2\pi}(u\ast P_r)$, u is in $C^p(\mathbb{T})$ and therefore $U\in C^p(\mathbb{T})$ and:
\begin{equation}
\frac{d^lU}{d\theta^l}(re^{i\theta_0})= \frac{1}{2\pi}(\frac{d^lu}{d\theta}\ast P_r)= \frac{1}{2\pi}\int_{-\pi}^{\pi}u^{(l)}(t)P_r(\theta_0- t) dt
\end{equation}
for $l= 0, 1, ..., p$, due to well known properties of convolution. The uniform convergence and (2.3) derives automatically from Poisson's kernel properties.

\end{proof}
\end{lemma}

\noindent Now we investigate for f in $A(D)$ what is the relation of the derivatives $f^{(l)}$ and the differentiability of the map $t\mapsto f(e^{it})$.
\begin{theorem}
For all $f\in A(D)$ the following equivalence holds: $f\in A^p(D)$ if and only if the map $g: \mathbb{R}\to \mathbb{C}$, defined by $g(t)= f(e^{it})$, is in $C^p(\mathbb{T})$. In that case:
\begin{equation}
\frac{dg}{dt}(t)= ie^{it}f'(e^{it})
\end{equation}

\begin{proof}

We will prove it by induction on $p$. For $p= 1$, let an $f\in A^1(D)$, define: $f_r(t)= f(re^{it}): \mathbb{R}\to \mathbb{C}$, for all $0<r<1$ and notice that $f_r\to g$ uniformly for $t\in\mathbb{R}$, since f is continuous in $\overline{D}$. Observe that:
\begin{equation}
\frac{df_r}{d\theta}(t)= ire^{it} f'(re^{it})= h_r(t)
\end{equation}
which uniformly converges to $h(t)= ie^{it}f'(e^{it})$ for $t\in\mathbb{R}$ as $r\to 1^-$ due to $f'$ being uniformly continuous in $\overline{D}$. Additionally, $f$ is continuous in $\overline{D}$ and therefore $f_r(0)= f(r)\to f(1)$, as $r\to 1^-$. From a well known proposition $f_r$ uniformly converges to $\int h(t)dt + c$ for a proper $c\in\mathbb{C}$ and hence $g(t) =\int h(t)dt + c$ which implies that $g\in C^1(\mathbb{R})$ and $g'(t)= h(t)= ie^{it}f'(e^{it})$.

For the converse let $f\in A(D)$ such that $g(t)= f(e^{it})\in C^1(\mathbb{R})$. It's well known that:
\begin{equation}
f(re^{i\theta})= \frac{1}{2\pi} \int_{-\pi}^{\pi}f(e^{i(\theta- t)})P_r(t)dt = \frac{1}{2\pi} \int_{-\pi}^{\pi}g(\theta- t)P_r(t)dt
\end{equation}
for all $0<r<1$ and $\theta\in\mathbb{R}$. 
By Proposition 2.1 we get:

\begin{equation}
ire^{i\theta}f'(re^{i\theta}) = \frac{df}{d\theta}(re^{i\theta})= \frac{1}{2\pi} \int_{-\pi}^{\pi} g'(t)P_r(\theta-t)dt \to g'(\theta)
\end{equation}
uniformly for $\theta\in\mathbb{R}$ and therefore:

\begin{equation}
\lim_{r\to 1^-}f'(re^{i\theta}) = \frac{g'(\theta)}{ie^{i\theta}}
\end{equation}
uniformly for $\theta\in\mathbb{R}$ and hence $f'$ extends continuously in $\overline{D}$ which by definition implies that $f\in A^1(D)$.
Now let us assume that the theorem holds for all $k< p$ for some $p\geq2$. For the straight direction, if $f\in A^p(D)$ then $f'\in A^{p-1}(D)$ and therefore $h(t)=f'(e^{it})\in C^{p-1}(\mathbb{T})$ by the induction hypothesis. Additionally, $f\in A^1(D)$ and hence $g\in C^1(\mathbb{T})$ with $g'(t)= ie^{it}f'(e^{it})\in C^{p-1}(\mathbb{T})$ which implies that $g\in C^p(\mathbb{T})$. For the converse, if $g\in C^p(\mathbb{T})$ then $g'(t)= ie^{it}f'(e^{it})\in C^{p-1}(\mathbb{T})$ which implies that the map $h(t)= f'(e^{it})\in C^{p-1 }(\mathbb{T})$. By the induction hypothesis we have that $f'\in A^{p-1}(D)$ and therefore $f\in A^p(D)$.
\end{proof}
\end{theorem}

\noindent As an interesting corollary we get the following:
\begin{corollary}
$A^p(D)= A(D) \cap  C^p(\mathbb{T})$.
\begin{proof}
The equality of the sets derives from proposition 2.2. Additionally, $C^p(\mathbb{T})$ is topologized through the semi-norms:
\[
|g|_l = \sup_{t\in\mathbb{R}}|g^{(l)}(t)|, 0\leq l\leq p, l\in\mathbb{N}
\]
and therefore, if $f\in A(D)$ and $g(t)= f(e^{it})$, it is not hard for one to see that for all $1\leq l\leq p, l\in\mathbb{N}$:
\begin{equation}
g^{(l)}(t)= \sum_{k=1}^{l}P_{k,l}(e^{it})f^{(k)}(e^{it})
\end{equation}
where $P_{k, l}(z)$ are fixed polynomials (see \cite{[2]} Lemma 2.15 for a similar result). It follows:
\begin{equation}
f^{(l)}(e^{it})= \sum_{k=1}^{l} Q_{k,l}(e^{-it})g^{(k)}(t)
\end{equation}
where $Q_{k, l}(z)$ are fixed polynomials. Let us denote by $a_{k,l}$ and $b_{k, l}$ the sum of the absolute values of the coefficients of $P_{k,l}(z)$ and $Q_{k, l}(z)$ respectively, for all $1\leq l\leq p, l\in\mathbb{N}$. Consequently, we have that:
\begin{equation}
|f|_0= \sup_{z\in\overline{D}}|f(z)|= \sup_{z\in\mathbb{T}}|f(z)|= \sup_{t\in\mathbb{R}}|f(e^{it})|= \sup_{t\in\mathbb{R}}|g(t)|= |g|_0
\end{equation}
and for $1\leq l\leq p, l\in\mathbb{N}$:
\begin{equation}
|f|_l= \sup_{z\in\overline{D}}|f^{(l)}(z)|= \sup_{z\in\mathbb{T}}|f^{(l)}(z)|\leq \sum_{k=1}^{l}b_{k, l}\sup_{t\in\mathbb{R}}|g^{(k)}(t)|= \sum_{k= 1}^{l}b_{k, l}|g|_k
\end{equation}
and:
\begin{equation}
|g|_l= \sup_{t\in\mathbb{R}}|g^{(l)}(t)|\leq \sum_{k=1}^{l}a_{k, l}\sup_{t\in\mathbb{R}}|f^{(k)}(e^{it})|= \sum_{k=1}^{l}a_{k, l}|f|_k
\end{equation}
and therefore the topologies induced are equivalent. 
\end{proof}

\end{corollary}

\section{Extendability over an open arc of the unit circle.}

The question that now naturally occurs is if functions in the $A(D)$ whose derivatives extend only on an open arc of the unit circle, but not necessarily on the whole circle, behave in a similar way. We will see that there is a positive answer to the question above. 
For that to be achieved, we will need to investigate the properties of functions that are defined as the convolution of integrable functions defined on the unit circle which are equal to $0$ everywhere but on the arc in which we are interested in, with the Poisson kernel. This method is natural to consider since the Poisson representation solved the problem for functions whose derivatives were defined on the whole unit circle.

\begin{proposition}
Let $u:[0, 1]\to \mathbb{R}$ a continuous function and define:
\begin{equation}
A(z) = \frac{1}{2\pi}\int_{0}^{1} u(t)P_z(t)dt 
\end{equation}
is well defined in $\mathbb{C}\setminus\{ e^{it}: 0\leq t \leq 1\}$ and $C^{\infty}$ harmonic.

\begin{proof}
We will first show that $A(z)$ is well defined on the above set. Observe that:
\begin{equation}
A(z) = \frac{1}{2\pi}\int_{0}^{1} u(t)P_z(t)dt = Re \left( \frac{1}{2\pi}\int_{0}^{1} u(t) \frac{1+e^{-it}z}{1-e^{-it}z} dt \right)
\end{equation}
and $1-e^{-it}z=0 \Leftrightarrow z=e^{i(t+2k\pi)}, k\in\mathbb{Z}$ and for a fixed $z\in\mathbb{C}\setminus\{ e^{it}: 0\leq t \leq 1\}$ observe that $\delta_z= dist(1, \{ze^{-it}: 0\leq t \leq 1\}>0$
and:
\begin{equation}
\left|  u(t) \frac{1+e^{-it}z}{1-e^{-it}z} \right| \leq \sup_{t\in[0, 1]}|u(t)| \frac{1+ |z|}{\delta_z} < +\infty
\end{equation} 
for all $t\in[0, 2\pi]$. Since the quantity on the right hand side is integrable we deduce that A(z) is indeed well defined in $\mathbb{C}\setminus\{ e^{it}: 0\leq t \leq 1\}$. 

It suffices to show that $g(z)= \frac{1}{2\pi}\int_{0}^{1} u(t) \frac{1+e^{-it}z}{1-e^{-it}z} dt, z\in\mathbb{C}\setminus\{ e^{it}: 0\leq t \leq 1\}$ is holomorphic and from (3.2) follows that $A(z)$ is $C^{\infty}$ harmonic. Notice that for $z\neq z_0$: 

\begin{equation}
\frac{g(z)- g(z_0)}{z- z_0}= \frac{1}{2\pi} \int_{0}^{1}u(t)\frac{2e^{-it}}{(1-e^{-it}z)(1-e^{-it}z_0)}dt
\end{equation}
and for $z$ sufficiently close to $z_0$: 

\begin{align}
&\left| \frac{g(z)- g(z_0)}{z- z_0} - \frac{1}{2\pi}\int_{0}^{1}u(t)\frac{2e^{-it}}{(1-e^{-it}z_0)^2}dt \right| = \\ \nonumber
&\left| \frac{1}{2\pi}\int_{0}^{1}u(t)\frac{2e^{-2it}(z-z_0)}{(1-e^{-it}z)(1-e^{-it}z_0)^2} dt\right|  \\ \nonumber
&\leq \frac{1}{2\pi} \sup_{t\in[0, 1]}|u(t)| 2 \left( \frac{2}{\delta_{z_0}} \right)^3 |z-z_0|  \xrightarrow{z\to z_0} 0
\end{align}
therefore g is holomorphic and the proof is complete.

\end{proof}
\end{proposition}

\begin{lemma}
Let $u:[0, 1]\to \mathbb{R}$ be a continuous function and A(z) as in Proposition 3.1. For all $z= re^{i\theta}\in \mathbb{C}\setminus\{ e^{it}: 0\leq t \leq 1\}$:
\begin{equation}
\frac{dA}{d\theta}({re^{i\theta}})= \frac{1}{2\pi} \int_{0}^{1}u(t)\frac{-2r(1-r^2)sin(\theta-t)}{(1+r^2- 2rcos(\theta-t))^2}dt
\end{equation}
and therefore for all $\theta\in (1, 2\pi)$ and $r=1$:
\begin{equation}
\frac{dA}{d\theta}({e^{i\theta}}) = 0
\end{equation}

\begin{proof}
Since $A(re^{i\theta}) = \frac{1}{2\pi}\int_{0}^{1} u(t)P_r(\theta- t)dt $ = $\frac{1}{2\pi} (f\ast P_r)$ and $P_r$ is differentiable in respect to $\theta$ we deduce that $A(re^{i\theta})$ is differentiable in respect to $\theta$ and:
\begin{align}
\frac{dA}{d\theta}\bigg|_{re^{i\theta}} &= \frac{1}{2\pi} \frac{d(f\ast P_r)}{d\theta}\bigg|_{re^{i\theta}}=
 \frac{1}{2\pi}(f\ast \frac{dP_r}{d\theta})\bigg|_{re^{i\theta}} \Longrightarrow \\ 
\frac{dA}{d\theta}\bigg|_{re^{i\theta}} &= \frac{1}{2\pi} \int_{0}^{1}u(t)\frac{-2r(1-r^2)sin(\theta-t)}{(1+r^2- 2rcos(\theta-t))^2}dt
\end{align}
and the proof is complete.
For (3.7) just apply (3.6) for $r=1$ and any $\theta\in (1, 2\pi)$.
\end{proof}
\end{lemma}

Since every continuous $u: [0, 1]\to\mathbb{C}$ can be considered as a $2\pi$- periodic function $v:\mathbb{R}\to\mathbb{C}$ such that $ v(x)= u(x)$ for $x\in [0, 1]+ 2\pi\mathbb{Z}$ and $v(x)= 0$ otherwise, one can expect that such a function when convolved with the Poisson kernel would retain the nice properties. More specifically, it will uniformly converge to $0$ and to $u(x)$ on the compact subsets of the respective open arcs. This is what we prove in Propositions 3.3 and 3.4.
\begin{proposition}
Let $A(z)$ defined as above, then for all $\theta\in (1, 2\pi)$:
\begin{equation}
\lim_{r\to 1^-}\frac{d^lA}{d\theta^l}(re^{i\theta})=0
\end{equation}
for all $l\in\mathbb{N}$, and the convergence is uniform in the compact subsets of $(1, 2\pi)$.
\begin{proof}
We will firstly prove it for $l= 1$. It suffices to prove that for all $[\theta_1, \theta_2]\subset (1, 2\pi)$ the convergence is uniform. Observe that $0<\theta_1 - t<\theta_2 -t< 2\pi$ for all $t\in[0,1]$ and therefore $cos(\theta- t)\leq \max\{cos(\theta_1 -1), cos(\theta_2) \}= M<1, \forall t\in[0,1]$ and $\forall \theta\in [\theta_1, \theta_2]$ which implies: $1+ r^2 -2rcos(\theta -t) \geq 1+ r^2 -2rM= (1- r)^2 +2r(1- M)$ for all $t\in [0, 1]$ and $\theta\in [\theta_1, \theta_2]$. We can now write for all $\theta\in [\theta_1, \theta_2]$ and $0< r< 1$:

\begin{align}
\left| \frac{dA}{d\theta}(re^{i\theta}) \right| = \left| \frac{1}{2\pi} \int_{0}^{1}u(t) \frac{-2r(1-r^2)sin(\theta-t)}{(1+r^2- 2rcos(\theta-t))^2}dt\right| \\
\leq \frac{1}{2\pi}  \norm{u}_{\infty} \frac{2r(1-r^2)}{(1- r)^2 +2r(1- M)} 
\end{align}
but (3.12) converges to $0$ as $r\to 1^-$ and hence $\sup\limits_{\theta\in [\theta_1, \theta_2]}\left| \frac{dA}{d\theta}(re^{i\theta}) \right| \xrightarrow{r\to 1^-} 0$. \\
\noindent Notice that no matter how many times we differentiate $P_r$ in respect to $\theta$ we will have a finite sum of fractions that consist of a power of $(1+ r^2- 2rcos(\theta- t))$ as a denominator and the product of a non zero real number, $1- r^2$ and powers of $r$, $cos(\theta- t)$ and $sin(\theta- t)$ as a numerator. As a result, for any $l\geq 2$ the same proof applies.
\end{proof}

\end{proposition}

\begin{proposition}
Let $p\in\{0, 1, ... \}\cup\{\infty\}$ and $u:(-\varepsilon, 1+ \varepsilon)\to \mathbb{R}$, for some $\varepsilon>0$, continuous such that $u\in C^p((-\varepsilon, 1+ \varepsilon))$ and define $A(z)$ as in proposition 3.1, then for all $\theta_0\in(0, 1)$:
\begin{equation}
\lim_{\substack{z\to e^{i\theta_0} \\ |z|<1}} \frac{d^lA}{d\theta^l}(z) = u^{(l)}(\theta_0)
\end{equation}
for $0\leq l\leq p, l\in\mathbb{N}$, while the convergence is uniform in the compact subsets of $(0, 1)$.

\begin{proof}
Let a $p\in \{0, 1, ...\}\cup\{\infty\}$. Due to Borel's theorem \cite{[5]} we can find a $q:\mathbb{R}\to\mathbb{R}$ function in $C^{\infty}(\mathbb{R})$ such that $q^l(0)= u^l(0)$ and $q^l(1)= u^l(1)$ for $0\leq l\leq p, l\in\mathbb{N}$. This implies that the function:
\[
g(t)=
\begin{cases}
u(t), t\in [0,1]\\ 
q(t), t\in(1, 2\pi)
\end{cases}
\]
is in $C^p(\mathbb{T})$ and define:
\begin{equation}
G(z)= \frac{1}{2\pi}\int_{0}^{2\pi} g(t)P_z(t)dt
\end{equation}
From Lemma 2.1 we know that uniformly for all $\theta\in \mathbb{R}$:
\begin{equation}
\lim_{\substack{z\to e^{i\theta} \\ |z|<1}} \frac{d^lG}{d\theta^l}(z) = g^{(l)}(\theta)
\end{equation}
for all $0\leq l\leq p, l\in\mathbb{N}$. Now we define:
\begin{equation}
A(z)= \frac{1}{2\pi}\int_{0}^{1} u(t)P_z(t)dt  \hspace{0.2cm} \text{and} \hspace{0.2cm}     B_l(z)= \frac{1}{2\pi}\int_{1}^{2\pi} q(t)P_z(t)dt
\end{equation}
notice that $G(z)= A(z)+ B(z)$ for all $|z|<1$ and hence for $0\leq l\leq p, l\in\mathbb{N}$:
\begin{equation}
\frac{d^lG}{d\theta^l}(z) = \frac{d^lA}{d\theta^l}(z) + \frac{d^lB}{d\theta^l}(z)
\end{equation}
Now let a $\theta_0\in (0,1)$, combining (3.15),  (3.17) and applying Proposition 3.3 for $B(z)$ we have that  $\lim\limits_{\substack{z\to e^{i\theta_0} \\ |z|<1}} \frac{d^lB}{d\theta^l}(z) = \frac{d^lB}{d\theta^l}(e^{i\theta_0})= 0 $  and:
\begin{equation}
\lim_{\substack{z\to e^{i\theta_0} \\ |z|<1}} \frac{d^lA}{d\theta^l}(z) = \lim_{\substack{z\to e^{i\theta_0} \\ |z|<1}}  \left( \frac{d^lG}{d\theta^l}(z) - \frac{d^lB}{d\theta^l}(z) \right) =
u^{(l)}(e^{i\theta_0})
\end{equation}
while the convergence is uniform on the compact subsets of $(0,1)$.
\end{proof}

\end{proposition}

\noindent Notice that Proposition 3.4 holds under the weaker assumption of a $u: [0, 1]\to \mathbb{R}$ such that $u\in C^p((0, 1))$ and all its derivatives $u^{(l)}$ for $0\leq l\leq p, l\in\mathbb{N}$, continuously extend on $[0, 1]$.

\begin{remark}
Propositions 3.1, 3.3, 3.4 and Lemma 3.2 hold for complex functions $f=u +iv$ where $u=Ref$ and $v= Imf$, since we can apply them for the Real and Imaginary part separately and due to the linearity of the integral.
\end{remark}

We are now in place to prove the first of the two main theorems of this section. Using the propositions and lemmas of this section we are able to mimic the proof of theorem 2.2 and eventually prove a similiar result for functions whose derivatives only extend over an open arc of  the unit circle.
\begin{theorem}
Let $p\in\{1, 2, ...\}\cup\{\infty\}$, $f\in A(D)$ and define $g: (0, 1) \to \mathbb{C}$, by setting $g(t)= f(e^{it})$, then  $f^{(l)}$ continuously extends on $D\cup \{ e^{i\theta}: 0<\theta< 1\}$ for all $0\leq l\leq p, l\in\mathbb{N}$ if and only if $g\in C^p((0,1))$. In that case, for all $t \in (0, 1)$:
\begin{equation}
\frac{dg}{dt}(t) = ie^{it}f'(e^{it})
\end{equation}

\begin{proof}
We will prove it by induction on $p$. For $p= 1$, let $t_0\in (0, 1)$ and $t_1<t_0<t_2$ such that $[t_1, t_2]\subset (0, 1)$. For the straight direction, define $f_r(t)= f(re^{it}): (0,1)\to \mathbb{C}$, $h(t)= ie^{it}f'(e^{it}): (0, 1)\to \mathbb{C}$ and:
\begin{equation}
h_r(t)= \frac{df_r}{dt}(t)= ire^{it}f'(re^{it})
\end{equation}
for all $t\in (0, 1)$, then $h_r\to h$ as $r\to 1^-$ uniformly in $[t_1, t_2]$ due to the continuity of $f'$ in $D\cup \{ e^{i\theta}: 0<\theta< 1\}$. Notice that $f_r(t_0)\to f(e^{it_0})$ and hence from a well-known theorem $f_r\to \int hdt + c$, for a proper $c\in\mathbb{C}$, as $r\to 1$ and the convergence is uniform on $[t_1, t_2]$. Additionally, $f_r\to g$, as $r\to 1$, uniformly in $[t_1, t_2]$ and therefore $g'(t)= h(t)= ie^{it}f'(e^{it}), t\in(t_1, t_2)$. Since $t_0\in (0, 1)$ was a random choice we deduce that $g\in C^1((0, 1))$ and its derivative is given by (3.19).

For the converse let $g\in C^1((0,1))$. We know that for all $|z|< 1$:
\begin{equation}
f(z)= \frac{1}{2\pi} \int_{0}^{2\pi} f(e^{it})P_z(t)dt
\end{equation}
Define:
\begin{equation}
A(z)= \frac{1}{2\pi} \int_{t_1}^{t_2} g(t)P_z(t)dt \text{ and }B(z)= \frac{1}{2\pi} \int_{t_2}^{2\pi + t_1} f(e^{it})P_z(t)dt
\end{equation}
then $f(z)= A(z)+ B(z)$ for all $|z|< 1$.  Therefore:

\begin{equation}
\frac{df}{dt}(re^{it})= \frac{dA}{dt}(re^{it}) + \frac{dB}{dt}(re^{it})
\end{equation}
for all $0<r<1$ and $t \in\mathbb{R}$. Since $[t_1, t_2]\subset (0, 1)$ Propositions 3.4 and 3.3 imply that:
\begin{equation}
\lim_{r\to 1^-} \frac{dA}{dt}(re^{it}) = \frac{dg}{dt}(t)
\end{equation}
 and 
 \begin{equation}
 \lim_{r\to 1^-} \frac{dB}{dt}(re^{it}) = 0.
 \end{equation}
uniformly for $t\in [t_1, t_2]$. Combining (3.23), (3.24) and (3.25) we get:
\begin{equation}
\lim_{r\to 1^-} f'(re^{it})= \lim_{r\to 1^-}\frac{1}{ire^{it}}\frac{df}{dt}(re^{it})  =\frac{1}{ie^{it}} \frac{dg}{dt}(t)
\end{equation}
uniformly for $t\in [t_1, t_2]$. Since $t_0$ was an arbitrary choice in $(0, 1)$ we deduce that $f'$ extends continuously on $D\cup \{e^{it}: 0< t< 1 \}$.
To complete the induction, let us assume that the theorem holds for some $p\geq 1$. For the straight direction, if $f^{(l)}$ continuously extends on $D\cup\{e^{it}: 0< t< 1\}$ for all $0\leq l\leq p+1, l\in\mathbb{N}$ it follows that $(f')^{(l)}$ continuously extends on $D\cup\{e^{it}: 0< t< 1\}$ for all $0\leq l\leq p, l\in\mathbb{N}$. By the induction hypothesis the map $t\mapsto f'(e^{it})$ belongs in the class $C^p((0, 1))$ and since, by the case of $p= 1$, $g'(t)= ie^{it} f'(e^{it})$ we have $g'\in C^p((0, 1))$ and hence $g\in C^{p+ 1}((0, 1))$. For the converse, if $g\in C^{p+ 1}((0, 1))$ it follows that $g'(t)= ie^{it}f'(e^{it})\in C^p((0, 1))$ and therefore the map $t\mapsto f'(e^{it})$ is in the class of $C^p((0, 1))$. By the induction hypothesis, $(f')^{(l)}$ continuously extends on $D\cup\{e^{it}: 0< t< 1\}$, for all  $0\leq l\leq p, l\in\mathbb{N}$ and hence $f^{(l)}$ continuously extends on $D\cup\{e^{it}: 0< t< 1\}$, for all $0\leq l\leq p+1, l\in\mathbb{N}$. The case of $p= \infty$ follows  easily.
\end{proof}
\end{theorem}

Using Riemann's mapping theorem \cite{[1]} and theorem 3.6 we are able to prove a more general version of the previous theorem. The function $f$ does not have to be in the class $A(D)$, but we can assume that $f$ continuously extends only on the open arc we are interested in.
\begin{theorem}
Let $p\in\{1, 2, ...\}\cup\{\infty\}$ and $f: D\cup \{e^{it}: 0<t <1 \} \to \mathbb{C}$ be continuous and holomorphic in $D$ and define $g: (0, 1)\to \mathbb{C}$ by setting $g(t)= f(e^{it})$. Then $f^{(l)}$ extends continuously on $D\cup \{e^{it}: 0<t <1 \}$, for all $0\leq l\leq p, l\in\mathbb{N}$, if and only if  $g\in C^p((0, 1))$. In that case:
\begin{equation}
\frac{dg}{dt}(t)= ie^{it}f'(e^{it})
\end{equation}
for all $t\in (0, 1)$.

\begin{proof}
We will prove it by induction on $p$. Let $p= 1$. For the straight direction the proof is no different than that of Theorem 3.6. From the straight direction we also obtain (3.27).  For the converse let a $t_0\in (0, 1)$ and $t_1< t_0< t_2$ such that $[t_1, t_2]\subset (0, 1)$ and consider: 
$V= \{|z|< 1: argz\in (t_1, t_2) \}$. It is easy to see that V is simply connected and hence there exists a Riemann mapping \cite{[1]} $\phi: D\to V$ which extends to homeomorphism over the closures $\overline{D}$ and $\overline{V}$, due to the Osgood- Caratheodory theorem \cite{[4]}. Since $\phi(\mathbb{T})= \{ e^{it}: t_1\leq t\leq t_2\}\cup \{re^{it_1}: 0\leq r\leq 1 \} \cup\{re^{it_2}: 0\leq r\leq 1  \}$ and $\{ e^{it}: t_1\leq t\leq t_2\}$ is connected we deduce that $\phi^{-1}(\{ e^{it}: t_1\leq t\leq t_2\})= \{e^{it}: 0\leq t \leq 1 \}$ where we picked $0, 1$ without any loss of generality. Therefore $f\circ\phi: \overline{D} \to \mathbb{C}$ is continuous on $\overline{D}$ and holomorphic in D. Using the Reflection Principle we deduce that $\phi$ has an analytic extension $\Phi$ in a neighborhood of $D\cup \{ e^{it}: t\in (0, 1) \}$, such that $\Phi'(z) \neq 0$ (see \cite{[1]} p. 233- 235) and therefore the map $t \mapsto \phi(e^{it})$ is in $C^{\infty}((0, 1))$ and its derivative is non zero. Since, by hypothesis, $g\in C^1((t_1, t_2))$ and the map $t \mapsto \phi(e^{it})$ is in $C^{\infty}((0, 1))$ this implies that the map $t\mapsto f(\phi(e^{it}))$ is in $C^1((0, 1))$. Due to theorem 3.6 the derivative $(f\circ\phi)'$ continuously extends on $D\cup \{e^{it}: 0< t< 1 \}$ and hence $f'$ continuously extends in $D\cup \{e^{it}: t_1< t< t_2 \}$, because $f= (f\circ\phi)\circ \phi^{-1}$ and $f'= (f\circ\phi)'\circ\phi^{-1}\cdot(\phi^{-1})'$. Notice that $t_0$ was an arbitrary choice and therefore $f'$ continuously extends in $D\cup \{e^{it}: 0< t< 1 \}$. To complete the induction the proof is no different than that of Theorem 3.6. Once again, the case of $p= \infty$ follows easily.
\end{proof}
\end{theorem}

\begin{remark}
Notice that 0 and 1, in the above theorems, propositions and lemmas can be replaced by any a, b respectively such that: $0\leq a< b < a+ 2\pi$, while the proofs remain exactly the same.
\end{remark}

\section{Jordan Domains bounded by an analytic curve.}
In this section we are interested in proving similar results over Jordan domains. We are able to prove analogous results on Jordan domains bounded by an analytic Jordan curve or a Jordan curve with an analytic part. We do this by making use of Riemann's mapping theorem combined with the results of the previous sections. Now let us give a definition:
\begin{definition}
Let $\gamma: I\to \mathbb{C}$ be a locally injective map  defined on the interval $I\subset\mathbb{R}$. Let also $t_0\in I^{\circ}$, where $I^{\circ}$ is the interior of $I$ in $\mathbb{R}$. A function $f: \gamma^{*}\to \mathbb{C}$ is one sided extendable at $(t_0, \gamma(t_0))$ if there exists a Jordan domain $\Omega$, such that $\partial\Omega$ contains an arc of $\gamma^{*}$, $\gamma[t_1, t_2]$ where  $t_1, t_2\in I$ such that $t_1< t_0< t_2$ and a continuous function $F:\Omega\cup\gamma(t_1, t_2)\to\mathbb{C}$ which is holomorphic in $\Omega$ and $F\restriction_{\gamma((t_1, t_2))}= f\restriction_{\gamma((t_1, t_2))}$.
\end{definition}

\begin{theorem}
Let $\gamma: I\to \mathbb{C}$ be a locally injective analytic curve, $t_0\in I^{\circ}$ and $f: \gamma^{*}\to\mathbb{C}$ a function which is one sided extendable at $(t_0, \gamma(t_0))$. Let $\Omega$ be the Jordan domain, $t_1, t_2\in I^{\circ}$ such that $t_1< t_0< t_2$ and $\gamma((t_1, t_2))\subset \partial\Omega$ and $F: \Omega\cup\gamma((t_1, t_2))\to\mathbb{C}$ the function which extends f according to definition 4.1. Define $g(t)= f(\gamma(t)): (t_1, t_2)\to \mathbb{C}$ and let $p\in \{1, 2,... \}\cup\{\infty\}$, then the following equivalence holds: $g \in C^p((t_1, t_2))$ if and only if $F^{(l)}$ continuously extends in $\Omega\cup\gamma((t_1, t_2))$ for all $0\leq l\leq p, l\in\mathbb{N}$. In this case we have:

\begin{equation}
\frac{dg}{dt}(t)= F'(\gamma(t)) \cdot \gamma'(t)
\end{equation}

\begin{proof}
We will prove it by induction on $p$. For $p=1$,  due to the fact that $\gamma$ is locally injective we can assume that $\gamma\restriction_{(t_1, t_2)}$ is one-to-one, for a proper choice of $t_1, t_2$. Since $\Omega$ is a Jordan domain there exists a Riemann map \cite{[1]}  $\phi: D\to \Omega$, which extends to a homeomorphism over the closures of $D$ and $\Omega$ according to the Osgood-Caratheodory theorem \cite{[4]}. Without any loss of generality we can assume that $\phi(\{e^{it}: 0< t< 1 \})= \gamma((t_1, t_2))$ and since $\gamma$ is an  analytic curve we know that the map $t \mapsto \phi(e^{it})$ is in $C^{\infty}((0, 1))$, according to the reflection principle $\cite{[1]}$. For the straight direction, let us assume that $g\in C^1((t_1, t_2))$. Notice that $h(t)= f(\phi(e^{it}))= (f\circ\gamma\circ\gamma^{-1}\circ\phi)(e^{it})= (g\circ \gamma^{-1}\circ \phi)(e^{it})$ and hence $h \in C^1((0, 1))$. Let $F$ denote the one-sided extension of $f$ onto $\Omega$. Since $F\circ\phi: D\cup \{e^{it}: 0<t <1\}\to \mathbb{C}$ is continuous on $\overline{D}$ and holomorphic in $D$ and additionally $H(t)= F(\phi(e^{it}))= h(t)\in C^1((0, 1))$ theorem 3.7 implies that $(F\circ\phi)'$ continuously extends on $D\cup\{e^{it}: 0< t<1 \}$ and:
\begin{align}
 \frac{dH}{dt}(t) &= ie^{it}(F\circ\phi)'(e^{it}) 
\end{align}
However $(F\circ\phi)'= (F'\circ\phi)\cdot\phi'$ for all $|z|< 1$ and since $\phi$ maps $\{e^{it}: 0< t< 1 \}$ to an analytic curve we know that $\phi$ has a holomorphic extension in a neighborhood $D\cup \{e^{it}: 0< t< 1\}$ such that $\phi'(z)\neq 0$ (see \cite{[1]} p. 233- 235), which implies that every derivative continuously extends on $D\cup\{e^{it}: 0< t< 1 \}$ and therefore $F'$ continuously extends in $\Omega\cup\gamma((t_1, t_2))$. Moreover, by theorem 3.7, for $t_0\in (t_1, t_2)$ we have:
\begin{equation}
\lim_{t\to t_0}\frac{(F\circ\phi)(\phi^{-1}(\gamma(t))- (F\circ\phi)(\phi^{-1}(\gamma(t_0))}{\phi^{-1}(\gamma(t))- \phi^{-1}(\gamma(t_0))} = (F\circ\phi)'(\phi^{-1}(\gamma(t_0)))
\end{equation}
and therefore:
\begin{gather}
\lim_{t\to t_0}\frac{f(\gamma(t))- f(\gamma(t_0))}{\gamma(t)- \gamma(t_0)}=  \\
\lim_{t\to t_0} \frac{(F\circ\phi)(\phi^{-1}(\gamma(t)))- (F\circ\phi)(\phi^{-1}(\gamma(t_0)))}{\phi^{-1}(\gamma(t))- \phi^{-1}(\gamma(t_0))} (\phi^{-1})'(\gamma(t_0))= \\
(F\circ\phi)'(\phi^{-1}(\gamma(t_0))) \cdot (\phi^{-1})'(\gamma(t_0))= \\
F'(\gamma(t_0)) \cdot \phi'(\phi^{-1}(\gamma(t_0))) \cdot (\phi^{-1})'(\gamma(t_0)) = F'(\gamma(t_0))
\end{gather}
Since $\gamma\in C^{\infty}((0, 1))$, it follows that:
\begin{equation}
\frac{dg}{dt}(t)= \frac{d(f\circ\gamma)}{dt}(t)= F'(\gamma(t))\gamma'(t)
\end{equation}

For the converse, since $F'$ continuously extends in $\Omega\cup\gamma((t_1, t_2))$ then we deduce that $F'\circ\phi$ continuously extends on $D\cup\{e^{it}: 0< t< 1 \}$ and therefore $(F\circ\phi)'= (F'\circ\phi)\cdot\phi'$ continuously extends on the same set. Once again, theorem 3.7 implies that $H(t)= F(\phi(e^{it}))= f(\phi(e^{it}))\in C^1((0, 1))$ and hence $g= f\circ\gamma= (f\circ\phi)\circ(\phi^{-1}\circ\gamma)\in C^1((t_1, t_2))$.\newline
To complete the induction let us assume that the theorem holds for all $k< p$ for some $p\geq 2$. For the straight direction, suppose that $g\in C^p((t_1, t_2))$ by the $k= 1$ case we have that $g'(t)= F'(\gamma(t))\cdot\gamma'(t)\in C^{p-1}((t_1, t_2))$ and since $\gamma'$ is continuous and $\gamma'(t)\neq 0$, we conclude that the map $t\mapsto F'(\gamma(t))$ is in $C^{p-1}((t_1, t_2))$. By the induction hypothesis $(F')^{(p-1)}= F^{(p)}$ continuously extends on $\Omega\cup\gamma((t_1, t_2))$. For the converse, if $F^{(l)}$ continuously extends on $\Omega\cup\gamma((t_1, t_2))$ for all $0\leq l\leq p, l\in\mathbb{N}$, then $(F')^{(l)}$ continuously extends on $\Omega\cup\gamma((t_1, t_2))$ for all $0\leq l\leq p-1, l\in\mathbb{N}$ and hence the induction hypothesis implies that $F'\circ\gamma\in C^{(p-1)}((t_1, t_2))$. But from the case $k= 1$ we deduce $\frac{d(F\circ\gamma)}{dt}= F'(\gamma(t))\cdot\gamma'(t)\in C^{p -1}((t_1, t_2))$ which implies that $F\circ\gamma\in C^p((t_1, t_2))$.This completes the proof for $p\neq\infty$. The case of $p=\infty$ follows easily.

\end{proof}
\end{theorem}
Let us now consider a Jordan domain $\Omega$ bounded by an analytic Jordan curve and let an $f\in A(\Omega)$. If $\gamma: \mathbb{R}\to\mathbb{C}$ is an analytic parametrization of $\partial\Omega$, we can consider the restriction of $f$ on $\partial\Omega$, $\tilde{f}: \gamma^\ast\to\mathbb{C}$. Notice that for all $(t, \gamma(t))$ then for proper $t_1< t< t_2$ we have $\gamma((t_1, t_2))\subset\partial\Omega$ and $f\in A(\Omega)$ such that $f\restriction_{(\gamma((t_1, t_2)))}$= $\tilde{f}\restriction_{(\gamma((t_1, t_2)))}$ and therefore $\tilde{f}$ is one-sided extendable at every point of $\gamma^\ast$ on the same Jordan domain $\Omega$ and by the same function $f\in A(\Omega)$. This in combination with theorem 4.2 implies that for an $f\in A(\Omega)$ the following equivalence holds: $f\in A^p(\Omega)$ if and only if $f\circ\gamma\in C^p(\mathbb{R})$. More specifically:

\begin{corollary}
Let $\Omega$ be a Jordan domain bounded by an analytic Jordan curve $\gamma: \mathbb{R}\to \partial\Omega$. For all $f\in A(\Omega)$ the following equivalence holds: $f\in A^p(\Omega)$ if and only if $f\circ\gamma\in C^p(\mathbb{R})$. In this case we have:
\begin{equation}
\frac{d(f\circ\gamma)}{dt}= f'(\gamma(t))\cdot\gamma'(t)
\end{equation}
\end{corollary}

If the Jordan domain $\Omega$ is bounded by an analytic Jordan curve and $\phi: D\to \Omega$ is a Riemann map, then it extends on a larger disk than $D$ \cite{[1]}. Then we use the parametrization $\gamma(t)= \phi(e^{it}), 0\leq t\leq 2\pi$. The topology of $C^p(\partial\Omega)$ is induced by the semi-norms
\begin{equation}
\sup_{0\leq t\leq 2\pi}\bigg|\frac{d^l(f\circ\gamma)}{dt^l}\bigg|, 0\leq l\leq p, l\in\mathbb{N}
\end{equation}
\noindent It follows easily that $A^p(\Omega)= A(\Omega)\cap C^p(\partial\Omega)$ and the proof is similar to the case $\Omega= D$.

\begin{remark}
The previous results can easily be extended in the case of finitely connected domains bounded by a finite set of disjoint Jordan curves, provided that a part of the boundary is analytic. Towards this a helpful tool is the extended Laurent Decomposition \cite{[3]}.
\end{remark}

\noindent We will take a moment to sketch the proof. Let $\Omega$ be a bounded domain whose boundary consists of a finite number of disjoint Jordan curves. If $V_0, V_1,... ,V_{n- 1}$ are the connected components of $\hat{\mathbb{C}}\setminus\Omega$, $\infty\in V_0$ and $\Omega_0= \hat{\mathbb{C}}\setminus V_0$, ..., $\Omega_{n-1}= \hat{\mathbb{C}}\setminus V_{n- 1}$, then for every $f$ which is holomorphic in $\Omega$ we know that there exist functions $f_0, ..., f_{n- 1}$ which are holomorphic in $\Omega_0, ..., \Omega_{n-1}$ respectively such that $f= f_0+ f_1+ ...+ f_{n- 1}$. Without any loss of generality, we assume that there exists an analytic arc on the boundary of $\Omega_0$ and denote by $\phi_0: \overline{D}\to \overline{\Omega_0}$ a Riemann map, holomorphic in $D$ and continuous on $\overline{D}$ \cite{[4]}. Since $f_1, ..., f_{n-1}$ are holomorphic in a neighborhood of $\partial\Omega_0$ we conclude that the derivatives of any order of the functions $f_1, ..., f_{n-1}$ continuously extend over $\partial\Omega_0$. Consequently, if $f$ continuously extends over an analytic arc of $\partial\Omega_0$ (equivalently $f_0$ continuously extends over the same arc) then from theorem 4.2 the first $p$ derivatives of $f$ (equivalently of $f_0$) continuously extend over that arc if and only if the map $t\mapsto (f\circ\phi_0)(e^{it})$ is in $C^p((t_1, t_2))$ for proper $t_1< t_2$ such that $\phi_0(\{e^{it}: t_1< t< t_2 \})$ is the analytic arc. 

\bigskip
\noindent \textbf{Acknowledgement}: The author would like to express his gratitude to
Professor Vassili Nestoridis for introducing him to the problem as well as for his guidance throughout the creation of this paper.

\bigskip

\noindent {\scshape Vlassis Mastrantonis} \\Department of Mathematics \\National and Kapodistrian University of Athens \\Panepistimioupolis, 157 84 \\Athens \\GREECE\\
E-mail: \href{mailto:blasismas@hotmail.com}{\url{blasismas@hotmail.com}}

\end{document}